\input amstex
\def\mapright#1{\smash{
   \mathop{\longrightarrow}\limits^{#1}}}

\def\mapdown#1{\Big\downarrow
   \rlap{$\vcenter{\hbox{$\scriptstyle#1$}}$}}

\def \cal{\Cal}
%AMS-TeX file "LaMu.tex"
\documentstyle{amsppt}
\NoRunningHeads \magnification=\magstep1 \baselineskip=12pt
\parskip=5pt
\parindent=18pt
\topskip=10pt \leftskip=0pt \rightskip=0pt \pagewidth{30pc}
\pageheight{47pc}

\topmatter
\title Low dimensional discriminant loci and scrolls
\endtitle
\author Antonio Lanteri and Roberto Mu\tildeaccent{n}oz \endauthor

\address Dipartimento di Matematica ``F. Enriques'',
Universit\`a, Via C. Saldini 50, I-20133 Milano, Italy
\endaddress
\email lanteri\@mat.unimi.it \endemail

\address
Departamento de Matem\'atica Aplicada, Universidad Rey Juan Carlos,
Calle Tulip\'an, 28933-M\'ostoles, Madrid, Spain
\endaddress
\email roberto.munoz\@urjc.es \endemail
%\date May 5, 2008
%\enddate

\abstract Smooth complex polarized varieties $(X,L)$ with a vector
subspace $V \subseteq H^0(X,L)$ spanning $L$ are classified under
the assumption that the locus ${\Cal D}(X,V)$ of singular elements
of $|V|$ has codimension equal to $\dim(X)-i$, $i=3,4,5$, the last
case under the additional assumption that $X$ has Picard number
one. In fact it is proven that this codimension cannot be
$\dim(X)-4$ while it is $\dim(X)-3$ if and only if $(X,L)$ is a
scroll over a smooth curve. When the codimension is $\dim(X)-5$
and the Picard number is one only the Pl\"ucker embedding of the
Grassmannian of lines in $\Bbb P^4$ or one of its hyperplane
sections appear. One of the main ingredients is the computation of
the top Chern class of the first jet bundle of scrolls and
hyperquadric fibrations. Further consequences of these
computations are also provided.
 \endabstract


%\subjclassyear{2000}
\endtopmatter

\document

\head Introduction
\endhead

Let $X$ be an irreducible smooth complex projective variety of
dimension $n$. Given a line bundle $L$ on $X$ and a linear system
$|V|$ of dimension $N$, defined by a vector subspace $V \subseteq
H^0(X,L)$, it is a classical problem to understand the (not
necessarily irreducible) subvariety ${\Cal D}(X,V)$ of $|V|$
parameterizing the singular elements of $|V|$. A first result on
this respect is Bertini theorem showing that, under the assumption
that $V$ spans $L$, the general element of $|V|$ cannot be
singular, so that ${\Cal D}(X,V)$ has positive codimension in
$|V|$. The hypothesis that $V$ spans $L$ seems then natural just
to have some control on the dimension of ${\Cal D}(X,V)$ (at least
to avoid the possibility of being the whole linear system).
Moreover, when $L$ is spanned by $V$ a morphism $\phi_V:X \to
\phi_V(X) \subseteq \Bbb P^N$ appears and the related geometry
enriches the picture. Imposing no other assumption on $V$, in
order to study the singular elements of $|V|$, we need to deal, on
the one hand, with the geometry of $\phi_V(X) \subset \Bbb P^N$
and, on the other, with the geometry of the fibres of $\phi_V$,
see for example \cite{LM1}. To overcome the problem represented by
positive dimensional fibres it is natural to impose ampleness on
$L$, so that $\phi_V$ is a finite morphism. According to
\cite{LPS1} these hypotheses, that is, $L$ ample with $V$ spanning
$L$, provide a good framework to study the object ${\Cal D}(X,V)$,
called the {\it discriminant locus} of the triplet $(X,L,V)$. In
fact, in this setting, ${\Cal D}(X,V)$ reflects a nice
stratification of $X$ determined by the rank of the differential
$d \phi_V$; in particular, a relevant role is played by the
geometry of $\phi_V(X) \subseteq \Bbb P^N$ and that of the
ramification locus of $\phi$ (for a more precise description see
Section 0: (0.1) and (0.2)).

If furthermore $L$ is very ample with $\phi_V$ defining an
embedding --from now on we refer to this context as the {\it
classical case}-- the discriminant variety is nothing else that
the dual variety $\phi_V(X)^\vee \subset \Bbb P^{N\vee}$, a very
classical object in projective geometry. For a survey on the topic
of dual varieties see \cite{T}.

A series of papers, see \cite{LPS1}, \cite{LPS2}, \cite{LM2},
\cite{LM3}, is devoted to investigate to what extent results holding
in the classical case are still true in the ample and spanned case.
Mostly, these works study the basic invariants, say dimension and
degree, of the discriminant locus. Let us focus on the dimension.
The hypothesis of $|V|$ being base point free allows us to write, by
Bertini theorem, $\dim({\Cal D}(X,V))=N-1-k$ where $k$ is called, as
usual, the {\it (discriminant) defect} of $(X,L,V)$. The problem of
classifying positive defect triplets appears naturally. In the
classical case a deep result by Beltrametti, Fania and Sommese, see
\cite{BFS, Thm. 1.2}, states that a positive defect variety is
always a fibration whose fibers are positive defect varieties with
Picard number one and such that their defect is even bigger than
that of the original variety. This shows that, when the defect is
very big with respect to the dimension, only scrolls can appear,
i.e., the fibers of the fibration are linear spaces. A concrete
result, always in the classical case, is that if $k \geq n-3>0$ then
$\phi_V(X)$ is either $\Bbb P^N$ or a scroll over a smooth curve
$C$. To avoid any confusion in the terminology, we say that $(X,L)$
is a scroll over a smooth variety $Y$ if there exists a vector
bundle $E$ over $Y$ such that $X=\Bbb P_Y(E)$ and $L$ is the
tautological line bundle.

In \cite{LPS1} and \cite{LM3} it is shown that, in the ample and
spanned case, the information on dimension and degree of the
discriminant variety is encoded in the Chern classes of the
so-called first jet bundle $J_1(L)$ of $L$. In particular, positive
defect is equivalent to the vanishing of the top Chern class of
$J_1(L)$. Let us consider the classification problem for triplets
$(X,L,V)$, with $L$ just ample and spanned by $V$, whose defect is
big with respect to $n$. Since scrolls $X=\Bbb P_Y(E)$ appear, at
least in the classical case, and by some evidences saying that they
seem to be the only examples, like \cite{LPS1, Conj. 2.11}, we need
to compute the top Chern class of $J_1(L)$ for scrolls. This is done
in Section 1, in fact without any assumption on $L$, and allows us
to compute in Section 2, see Proposition 2.1, the defect of scrolls
in the range $n-2\dim(Y) \geq -1$. Also irreducibility of the
discriminant locus is shown when $n-2\dim(Y)\geq 0$ and a
description of $\cal D$ when $(X,L)$ is the conormal variety of a
smooth curve is provided. In Section 3, see Theorems 3.1 and 3.2, we
prove that if the defect is bigger than or equal to $n-3$ (and
positive) then $(X,L)$ either is a scroll over a curve or
$(X,L)=(\Bbb P^n, \Cal O_{\Bbb P^n}(1))$. Moreover, if the defect is
equal to $n-4>0$ and the Picard number of $X$ is one then $(X,L)$ is
either $(\Bbb G(1,4),L)$, where $\Bbb G(1,4)$ is the Grassmannian of
lines in $\Bbb P^4$ and $L$ defines the Pl\"ucker embedding $\Bbb
G(1,4) \subset \Bbb P^9$, or there exists $H \in |L|$ such that
$(X,L)=(H, L_H)$, where $L_H$ stands for the restriction of $L$ to
$H$.

Our results of Section 3 rely on two basic facts. First one is the
existence of linear systems $|W| \subset |V|$ not meeting the
discriminant, that is, not containing any singular element. These
systems produce exact sequences involving $J_1(L)$ that lead to
the computation of its Chern classes. Second one is adjunction
theory, more specifically, the classification of polarized
varieties $(X,L)$ such that $K_X+(n-i)L$ is not ample for
$i=0,1,2$, see \cite{I} and \cite{F2}. For instance, the use of
adjunction makes hyperquadric fibrations over a curve enter into
the picture. This is the reason why we need to compute the top
Chern class of $J_1(L)$ also for such polarized varieties. This
computation is done in Section 1 and applied to exclude the
possibility that defect is $n-3$ in Section 3. Here we embed the
hyperquadric fibration into a scroll, fiberwise, as a divisor.

\noindent \eightpoint {\bf Acknowledgments.} We would like to thank
Raquel Mallavibarrena for several helpful discussions and the
referee for useful remarks. The first author is partially supported
by MiUR of the Italian Government (Cofin 2006, in the framework of
the National Research Project "Algebraic varieties, motives, and
arithmetic geometry") and by the University of Milano (FIRST 2005).
The second author is partially supported by the project Proyecto
MTM2006-04785 of the Spanish Government. \tenpoint

\head 0. Background material
\endhead

\flushpar {\bf (0.0)} Let $(X,L,V)$ be a triplet where: $X$ is an
irreducible smooth complex projective variety of dimension $n$,
$L$ is an ample and spanned line bundle on $X$ and $V\subseteq
H^0(X,L)$ spans $L$. Write $\dim(V)=N+1$ and let $\phi_V:X \to
\Bbb P^N$ (respectively $\phi_L$) be the morphism defined by $V$
(respectively by $H^0(X,L)$).

We use standard notation in algebraic geometry and the base field
is always $\Bbb C$. In particular, $K_X$ will denote the canonical
bundle of $X$.

\flushpar{\bf (0.1)} The {\it{discriminant locus}} $\Cal D(X,V)$
of the triplet $(X,L,V)$ parameterizes the singular elements of
$|V|$. More precisely, taking the incidence correspondence
$$\matrix
\Cal Y:=\{(x,[s]) \in X \times |V|:j_1(s)(x)=0\} & \mapright{p_1}
& X  \cr \mapdown{p_2} \cr \phantom{,,,,,,}\Cal D(X,V) \subset
\Bbb P^{N \vee},
\endmatrix$$
where $j_1(s)$ denotes the first jet of the section $s \in V$, $\Cal
D(X,V)$ is the image of $\Cal Y$ via the second projection of $X
\times |V|$. Thus $\Cal D(X,V)$ is an algebraic variety in $|V|=\Bbb
P^{N \vee}$. By Bertini theorem, $\dim (\Cal D(X,V)) < N$. Hence, we
can write $\dim (\Cal D(X,V)) = N-1-k$ where def$(X,V):=k \geq 0$ is
called the {\it{(discriminant) defect}} of $(X,L,V)$.  In particular
def$(X,V)$ is the dimension of the largest linear projective space
of $|V|$ not meeting ${\Cal D}(X,V)$.

If $\phi_V(X) \not= \Bbb P^N$ then the dual variety
$\phi_V(X)^{\vee} \subset \Bbb P^{N\vee}$ is a non-empty
irreducible subvariety of $\Cal D(X,V)$. Furthermore, if $\phi_V$
is an immersion then $\phi_V(X)^{\vee} = \Cal D(X,V)$, see
\cite{LPS1, Rmk. 2.3.3}. Anyway, points in $\Cal D(X,V) \setminus
\phi_V(X)^{\vee}$ are coming from points on $X$ where the
differential of $\phi_V$ is not injective. In this context it is
natural to define the {\it jumping sets} $\Cal J_i$  ($1 \leq i
\leq n$) as in \cite{LPS1, (1.1)}:
$$\Cal J_i(X,V)=\{x \in X: \text{rk}(d\phi_V(x)) \leq n-i\}.$$
As in \cite{LPS2, (0.3.1)} $X_i$ stands for $\Cal J_i \setminus
\Cal J_{i+1}$, with the convention that $\Cal J_0=X$ and $\Cal
J_{n+1}=\emptyset$. This allows to define $\Cal D_i(X,V) \subseteq
\Cal D(X,V)$ as $\overline{p_2\circ p_1^{-1}(X_i)}$, that is, the
Zariski closure in $\Bbb P^N$ of the locus of elements of $|V|$
singular at points of $X_i$. This gives a sort of stratification
of the discriminant locus, say
$$ \Cal D(X,V) = \cup_{i=0}^n \Cal D_i(X,V), \tag 0.2$$
and we can define def$_i(X,V)=N-1-\dim(\Cal D_i(X,V))$.

Let us observe that, as a consequence of the Bertini theorem, for
any choice of a vector subspace $W \subset V$ spanning $L$ we get
$\text{def}(X,V)=\text{def}(X,W)$, see \cite{LPS1, Thm. 2.7}. Hence,
as suggested after (2.7.2) in \cite{LPS1}, we can write def$(X,L)$
instead of def$(X,V)$, this value being independent of the choice of
$V$ spanning $L$. On the contrary the jumping sets are tightly
related to the choice of the linear system $|V|$, think, for
example, about projections of embedded projective varieties.

As shown in the definition of the discriminant locus and as will be
shown in (0.4) the first jet bundle $J_1(L)$ is a very remarkable
object in the study of the dimension and degree of the discriminant
locus. For the sake of completeness let us recall briefly its
construction, see \cite{LPS1, (0.3)} and reference therein. Let $x
\in X$ and let $\frak{m}_x$ be the maximal ideal of the stalk $\Cal
O_{X,x}$ at $x$ of the structure sheaf of $X$. Then
$\frak{m}_x^2L_x$ is the $\Cal O_{X,x}$ submodule of $L_x$ of the
germs of sections of $L$ whose local expansion at $x$ has no terms
of degree less than or equal to one. The quotient
$J_1(L)_x=L_x/\frak{m}_x^2L_X$ is a $\Bbb C$-vector space of
dimension $n+1$ and $J_1(L)$ is the collection $\cup_{x \in
X}J_1(L)_x$ equipped with the unique holomorphic vector bundle
structure inducing on $J_1(L)_x$ the preexistent $\Bbb C$-vector
space structure. Locally near $x$ a section of $J_1(L)$ is given by
a pair $(s,ds)$ where $s$ is a local section of $L$ and $ds$ is the
differential of $s$ with respect to local coordinates of $X$ around
$x$. The map sending $(s,ds)$ to $s$ gives rise to the following
exact sequence, where $\Omega_X$ is the contangent bundle:
$$0 \to \Omega_X (L) \to J_1(L) \to L  \to 0. \tag{0.3}$$

\noindent{\bf (0.4)} The dimension and the degree of the
discriminant locus can be computed by means of the Chern classes of
the first jet bundle. In fact, let us recall that def$(X,L) \geq r$
if and only if $c_{n-r+1}(J_1(L))=0$, see \cite{LPS1, Thm. 2.7}. In
particular, for $(X,L)$ to be defective the top Chern class
$c_n(J_1(L))$ must vanish. On the other hand, when $(X,L)$ is not
defective, $c_n(J_1(L))$ can be interpreted in terms of the degree
of the discriminant locus, usually called codegree and denoted by
codeg$(X,L)$. For this interpretation see \cite{LM3, Thm. 5.2}. Just
recall here that if the general section in any maximal dimensional
component of $\Cal D(X,L)$ is singular at a single point and the
singularity is ordinary quadratic of maximal rank then
$c_n(J_1(L))=$ codeg$(X,L)$. When equality $c_n(J_1(L))=$
codeg$(X,L)$ holds, for example when $L$ is very ample, we say that
$(X,L)$ has {\it tame codegree}, see \cite{LM3, Def. 10.1}.

\head 1. Top Chern classes of jet bundles
\endhead

As said in (0.4), to understand the dimension and degree of the
discriminant locus of a triplet $(X,L,V)$ as in (0.0) we need to
compute the top Chern class $c_n(J_1(L))$ of the first jet bundle.
This section is devoted to compute this Chern class for scrolls and
for hyperquadric fibrations over curves. In the case of scrolls no
assumption on the tautological line bundle is needed to do the
computation. However, in Section 2 we will assume ampleness and
spannedness to get nice results on the dimension of the discriminant
locus of scrolls. In the case of hyperquadric fibrations we compute
the top Chern class of the first jet bundle embedding our variety as
a divisor into a scroll.

First we compute $c_n(J_1(L))$ for scrolls. Recall the usual
convention that ${m \choose n}=0$ if either $n<0$ or $n>m$.

\proclaim{Proposition 1.1} Let $Y$ be an irreducible smooth
projective variety of dimension $m$ and $E$ a vector bundle of
rank $r$ on $Y$. Consider the projective bundle $\pi:\Bbb P_Y(E)
\to Y$ and the corresponding tautological line bundle $L$. Let
$n=m+r-1$ be the dimension of $X:=\Bbb P_Y(E)$. In these
conditions:
$$c_{n}(J_1(L))=\sum_{s_1,s_2 \geq 0,\; s_1+s_2 \leq m}A_{s_1,s_2}L^{n-s_1-s_2}
\pi^*c_{s_1}(E^\vee)\pi^*c_{s_2}(T_Y)$$ where
$$A_{s_1,s_2}=\sum_{t=0}^{n}(-1)^t(n+1-t){r-s_1 \choose
t-s_1-s_2}.$$
\endproclaim
\demo{Proof} By (0.3), see \cite{BS, Lemma 1.6.4}, we know that
$$c_n(J_1(L))=\sum_{t=0}^n(n+1-t)(-1)^tc_t(T_X)L^{n-t}. \tag{1.1.1}$$
Consider the relative Euler and  tangent exact sequences:
$$0 \to {\Cal O}_X \to \pi^*(E^\vee) \otimes L \to T_{X/Y} \to
0, \tag{1.1.2}$$ $$0 \to T_{X/Y} \to T_X \to \pi^*T_Y \to 0.$$ By
(1.1.1) and (1.1.2) we get:
$$c_{n}(J_1(L))=\sum_{t=0}^n(-1)^t(n+1-t)L^{n-t}
(\sum_{i=0}^tc_i(\pi^*(E^\vee) \otimes L)c_{t-i}(\pi^*(T_Y))).
\tag{1.1.3}$$ We end the proof computing the coefficient of
$L^{n-s_1-s_2}\pi^*c_{s_1}(E^\vee)\pi^*c_{s_2}(T_Y)$ in
$c_n(J_1(L))$ after substituting $c_i(\pi^*(E^\vee)\otimes L)$ in
(1.1.3) by
$$c_i(\pi^*(E^\vee)\otimes L)=\sum_{j=0}^i {r-j \choose
i-j}c_j(\pi^*(E^\vee))L^{i-j}. \qquad {_\square}$$  \enddemo

We can derive some consequences of Proposition 1.1 that will be
useful in Section 2.

\proclaim{Proposition 1.2} In the conditions of $(1.1)$ we get:

$(1.2.1)$ if $r \geq m+2$ then $c_n(J_1(L))=0$,

$(1.2.2)$ if $r=m+1$ then $c_n(J_1(L))=c_m(E)$,

$(1.2.3)$ if $r=m$ then $c_n(J_1(L))=c_{m-1}(E)(c_1(E)+K_Y)+m
c_m(E)$.
\endproclaim

\demo{Proof} Consider $f:\Bbb Z^4 \to \Bbb Z$ defined as
$$f(m,r,s_1,s_2)=\sum_{t=0}^{m+r-1}(-1)^t(m+r-t){r-s_1 \choose
t-s_1-s_2}.$$ We claim that:

(1.2.4) if $m,r,s_1,s_2 \geq 0$, $s_1+s_2 \leq m$ and $r \geq m+2$
then $f(m,r,s_1,s_2)=0$,

(1.2.5) if $m,r,s_1,s_2 \geq 0$, $s_1+s_2 \leq m$ and $r =m+1$
then $f(m,m+1,s_1,s_2)=0$ unless $s_1=m$ and $s_2=0$, in which
case $f(m,m+1,m,0)=(-1)^m$,

(1.2.6) if $m,r,s_1,s_2 \geq 0$, $s_1+s_2 \leq m$ and $r =m$ then
$f(m,m,s_1,s_2)=0$ unless $(s_1,s_2)
\in\{(m,0),(m-1,1),(m-1,0)\}$, in which cases $f(m,m,m,0)=(-1)^m
m$, $f(m,m,m-1,1)=(-1)^m$, $f(m,m,m-1,0)=(-1)^{m+1}$.

Let us observe that (1.2.1) follows directly from (1.2.4) and
(1.1). From (1.2.5) and (1.1) we get
$c_n(J_1(L))=(-1)^mL^{n-m}c_m(\pi^*(E^\vee))=c_m(E)$ so that
(1.2.2) holds. From (1.2.6), (1.1) and basic properties of Chern
classes we get:
$$c_n(J_1(L))=L^m\pi^*c_{m-1}(E)-L^{m-1}\pi^*c_{m-1}(E)\pi^*c_1(T_Y)+mL^{m-1}\pi^*c_m(E).$$
Finally use the Chern--Wu relation
$\sum_{i=0}^m(-1)^{i}\pi^*c_i(E)L^{m-i}=0$ to substitute $L^m$ in
the previous formula and get (1.2.3).

Let us now prove (1.2.4), (1.2.5) and (1.2.6). The proof is based
on two basic formulae, holding for $m \geq 2$:
$$\sum_{t=0}^{m-1} (-1)^t {m-1 \choose t}=\sum_{t=0}^m(-1)^t t {m \choose
t}=0. \tag{1.2.7}$$ First consider the decomposition of
$f(m,r,s_1,s_2)$ as a difference:
$$\sum_{t=0}^{m+r-1}(-1)^t(m+r){r-s_1 \choose
t-s_1-s_2}-\sum_{t=0}^{m+r-1}(-1)^t t{r-s_1 \choose t-s_1-s_2}.$$
Now make the change of variables $z=t-s_1-s_2$. If $(*)\;  s_2 \ne
m $ and $(**) \; r-s_1 \geq 2$ then by (1.2.7)
$$0=(-1)^{s_1+s_2}f(m,r,s_1,s_2)=$$ $$(m+r)\sum_{z=0}^{r-s_1}(-1)^z
{r-s_1 \choose z}-\sum_{z=0}^{r-s_1}(-1)^z (z+s_1+s_2){r-s_1
\choose z}.$$  If $r \geq m+2$ then (**) holds and (*) holds
unless $(s_1,s_2)=(0,m)$. Thus, to obtain (1.2.4) we only need to
check that $f(m,r,0,m)=0$, which is an easy computation. If
$r=m+1$ then (**) holds unless $s_1=m$. Thus, we obtain (1.2.5)
just checking that $f(m,m+1,m,0)=(-1)^m$ and $f(m,m+1,0,m)=0$.
Finally, if $r=m$ then (**) holds unless $s_1=m$ or $s_1=m-1$. We
get (1.2.6) just checking that $f(m,m,m,0)=(-1)^m m$,
$f(m,m,m-1,1)=(-1)^m$, $f(m,m,m-1,0)=(-1)^{m+1}$ and
$f(m,m,0,m)=0$. \qed \enddemo

Let us do a couple of remarks on this proposition. First we note
that (1.2.3) is related with the invariant $v(Y,E)$ defined by
Fukuma, see \cite{Fu, Section 4}:
$$v(Y,E)=1+{1\over 2}((m-2)c_m(E)+(K_Y+c_1(E))c_{m-1}(E)).$$
By (1.2.3) we get the following formula relating
the top Chern classes of $J_1(L)$ and $E$ with the invariant
$v(Y,E)$:
$$c_n(J_1(L))=2v(Y,E)-2+2c_m(E), \tag{1.2.8}$$ which will be of
interest in (2.1.2).

Second we note that the computation of $c_n(J_1(L))$ can be done for
$r<m$. However, formulae are more involved and difficult to control.
For example, if $r=m-1$ there are six non-zero coefficients in the
expression of $c_n(J_1(L))$. Let us recall that $c_n(J_1(L)) \geq
0$, see \cite{LPS1, Cor 2.6}, for $L$ ample and spanned. Hence,
these formulae provide non-negative expressions involving the Chern
classes of $E$ and $T_Y$.

This concludes our computations on scrolls. Let us start the
computation for hyperquadric fibrations over curves. We can use the
following notion of hyperquadric fibration over a smooth curve, see
\cite{F1, \S 3}.

\proclaim{Definition} Let $X$ be a smooth projective variety of
dimension $n$ and $L$ an ample line bundle on $X$. The pair $(X,L)$
is a {\rm hyperquadric fibration over a smooth curve} $B$ if there
exists a morphism $f:X \to B$  such that any general fiber $F$ of
$f$ is a smooth hyperquadric of $\Bbb P^n$ and $L$ induces the
hyperplane bundle on $F$, i.e., $L_F={\Cal O}_F(1)$.
\endproclaim

As in \cite{F1, (3.2)} it is a natural construction to take the
push-forward $E=f_*L$ which is a vector bundle of rank $n+1$ onto
$B$. This allows us to regard $X$ as a divisor in the projective
bundle $\pi:\Bbb P_B(E) \to B$. In fact, denoting by $\xi$ the
tautological bundle on $\Bbb P_B(E)$, there exists a divisor
$\beta \in $Pic$(B)$ of degree $b$ such that $X \in
|2\xi-\pi^*\beta|$ and $L = \xi_X$. Set $e=c_1(E)$ and recall that
$L^n=2e-b$.

\proclaim{Proposition 1.3} Let $(X,L)$ be a hyperquadric fibration
$f:X \to B$ over a smooth curve $B$ of genus $g$. Let $n=\dim(X)$
and let $e$ and $b$ be as before. Then $c_n(J_1(L))=2Ae-Bb
+4C(1-g)$, being
$$A=\sum_{t=0}^{n}(-1)^t(n+1-t)\big(\sum_{i=0}^t(-1)^{i}2^{i}{n
\choose t-i}\big), $$
$$B=\sum_{t=0}^{n}(-1)^t(n+1-t)\big(\sum_{i=0}^t(-1)^{i}2^{i}(i+1){n+1 \choose t-i}\big),$$
$$C=\sum_{t=0}^{n}(-1)^t(n+1-t)\big(\sum_{i=0}^t(-1)^{i}2^{i}{n+1
\choose t-i-1}\big).$$
\endproclaim

\demo{Proof} For simplicity write $P=\Bbb P_B(E)$ where $E=f_*L$.
The exact sequence associated to the embedding of $X$ as a divisor
in $P$ is $0 \to T_X \to (T_{P})_X \to \cal O_X(2\xi-\pi^*\beta)
\to 0$ and produces the following recursion law
$c_t(T_X)=c_{t}(T_P)_X-c_{t-1}(T_X)(2\xi-\pi^*\beta)_X$, for
$t=1,\dots,n$. This leads to the formula
$$c_t(T_X)=\sum_{i=0}^t(-1)^{i}c_{t-i}(T_P)_X(2\xi-\pi^* \beta)_X^{i}. \tag{1.3.1}$$
Now observe that
$(2\xi-\pi^*\beta)_X^{i}=2^{i}L^{i}-ib2^{i-1}L^{i-1}F$ and use the
exact sequences on $P$ analogous to (1.1.2)  to get
$$c_s(T_P)_X=(c_{s-1}(\pi^*(E^\vee) \otimes
\xi)c_1(\pi^*(T_B)))_X+c_s(\pi^*(E^\vee) \otimes \xi)_X,$$ for
$s=1, \dots, n$. Since
$$c_s(\pi^*(E^\vee)\otimes\xi)_X={n+1 \choose s}L^s-e{n \choose
s-1}FL^{s-1}$$ we can substitute in $c_s(T_P)_X$ to get
$$c_s(T_P)_X={n+1 \choose s}L^s-e{n \choose
s-1}FL^{s-1}+2(1-g){n+1 \choose s-1}L^{s-1}F.$$ Then (1.3.1) gives:
$$c_t(T_X)=\sum_{i=0}^t\big[ {n+1 \choose t-i}2^{i}L^{t}-{n+1 \choose
t-i}i2^{i-1}b L^{t-1}F-$$ $$ e{n \choose t-i-1}2^{i}L^{t-1}F
+2^{i+1}(1-g){n+1 \choose t-i-1}L^{t-1}F\big].$$ Finally we use
(1.1.1) together with $L^{n-1}F=2$ and $L^n=2e-b$ to get the
result. \qed
\enddemo

As an example, let us write the formula for $n=4$. We will use it
in (3.2.3). Let $(X,L)$ be as in (1.3). If $\dim(X)=4$ then, with
the same notation as there, $A=4$, $B=16$ and $C=-1$, so that
$c_4(J_1(L))=8e-16b-4(1-g)$. In fact we can rewrite it in the
following equivalent way
$$c_4(J_1(L))=(2e-b)+3(2e-5b)-4(1-g), \tag{1.4}$$ where
$2e-b =L^4>0$ by the ampleness of $L$ and $2e-5b\geq 0$ because it
corresponds to the number of singular fibers of $f$, see \cite{F1,
(3.3)}.

\head 2. Discriminant loci of scrolls
\endhead

Let us study the discriminant loci of scrolls. Recall the notation
established in the introduction: a scroll $(X,L)$ is a projective
bundle $\pi: X \to Y$ polarized by its tautological line bundle,
i.e., $X=\Bbb P_Y(E)$ where $E$ is an ample vector bundle of rank
$r$ on $Y$ and $L$ is the tautological line bundle. It is well known
in the classical case (i.e., when $L$ is very ample, see for
example, \cite{T, Thm. 7.21, p.\ 129}) that the defect of a scroll
of dimension $n$ over a $m$-dimensional base is greater than or
equal to $n-2m$ and in fact equal when $n-2m \geq 0$. Using the
computations of Proposition 1.2, since the vanishing of
$c_n(J_1(L))$ is needed to have positive defect, we give a proof of
this fact which works in the ample and spanned setup. We can also
deal with the case $n-2m=-1$ using the reformulation of (1.2.3)
stated in (1.2.8) combined with a result of Wisniewski \cite{Wi,
Thm. 3.4} classifying ample and spanned vector bundles with
$c_m(E)=1$. On the other hand, a result on the rank of this type of
vector bundles (apparently implicit in Wisniewski's result) can be
proved by our previous results as we show in the following remark.

\proclaim{Remark 2.0} {\rm Let $Y$ be an irreducible smooth
projective variety of dimension $m$ and $E$ a rank $r$ vector bundle
on $Y$. If $E$ is ample and spanned and $c_m(E)=1$ then $r=m$. In
fact, if $r<m$ then $c_m(E)=0$. If $r>m+1$ then we can choose
$r-(m+1)$ general independent sections which give rise a
homomorphism $\Cal O_Y^{\oplus r-(m+1)} \to E$ of maximal rank in
any fibre. This produces an exact sequence $$0 \to \Cal O_Y^{\oplus
r-(m+1)} \to E \to Q \to 0,$$ where $Q$ is a vector bundle which is
ample and spanned, so being $E$. Consider the pair $(X=\Bbb
P_Y(Q),L)$ where $L$ is the tautological line bundle and write
$n=\dim(X)$. By (1.2.2) we get $c_n(J_1(L))=c_m(E)=1$, that is, see
\cite{LM3, Thm. 5.2}, codeg$(X,L)=1$, contradicting \cite{LM3, Thm.
5.13}.}
\endproclaim

Let us introduce a bit of notation: Let $Y$ be an irreducible smooth
projective variety of dimension $m$, $E$ an ample vector bundle of
rank $r>1$ on $Y$ and $X=\Bbb P_Y(E) \mapright{\pi} Y$, $\dim(X)=n$.
Suppose that the tautological line bundle $L$ is spanned by $V
\subseteq H^0(X,L)$, $\dim(V)=N+1$. Consider the following incidence
variety $I \subset |V| \times Y$:
$$\matrix I & =\{(H,y):\; \pi^{-1}(y) \subset H\} & \mapright{p_1}
& |V|=\Bbb P^N \cr \mapdown{p_2} \cr \phantom{.}Y. \endmatrix \tag
2.1.0$$

\proclaim{Proposition 2.1} With the notation of $(2.1.0)$ we get
that:

$(2.1.1)$ If $n-2m \geq 0$ then {\rm def}$(X,L) =n-2m$, $p_1$ is
generically finite and $\Cal D(X,V)=p_1(I)$.

$(2.1.2)$ If $n-2m=-1$ and $m \geq 3$ then {\rm def}$(X,L) =0$
unless $Y=\Bbb P^m$ and $E=\cal O_{\Bbb P^m}(1)^{\oplus m}$, in
which case $(X,L) = (\Bbb P^m \times \Bbb P^{m-1}, {\cal O}_{\Bbb
P^m \times \Bbb P^{m-1}}(1,1))$.
\endproclaim

\demo{Proof} Let us first show the inequality {\rm def}$(X,L) \geq
n-2m$. Since def$(X,L) \geq 0$, we can suppose $n>2m$. By
\cite{LPS1, Thm. 2.7} it is enough to prove that
$c_{2m+1}(J_1(L))=0$. Consider a general element $X_1 \in |L|$ and
denote $L_1=L|_{X_1}$. In \cite{LM2, Lemma 1.13} it is shown that
$c_{2m+1}(J_1(L_1))=c_{2m+1}(J_1(L))|_{X_1}$ if $n>2m+1$. Hence we
can suppose $n=2m+1$, just taking general elements of $|L|$. We thus
conclude by (1.2.1).

If the strict inequality def$(X,L) > n-2m \geq 0$ holds then
$c_{2m}(J_1(L))=0$, again by \cite{LPS1, Thm. 2.7}. Exactly as in
the previous paragraph we can suppose $n=2m$ (just taking general
elements of $|L|$), i.e., $r=m+1$. By (1.2.2) we get
$0=c_{2m}(J_1(L))=c_m(E)$, which contradicts the ampleness of $E$,
see \cite{BG, Cor. 1.2}. This shows the first part of (2.1.1).

Now we observe that the incidence variety $I$ of (2.1.0) is a
projective bundle over $Y$ whose fiber is a linear space of
dimension $N-1-(n-m)$, in particular it is irreducible and of
dimension $N-1-(n-2m)$. Hence $p_1(I)$ is irreducible and
$\dim(p_1(I)) \leq N-1-(n-2m)$. Furthermore, since $r>1$, $\Cal
D(X,V) \subseteq p_1(I)$ so that $\dim(\Cal D(X,V)) \leq
N-1-(n-2m)$ or, equivalently, def$(X,L) \geq n-2m$. Moreover, the
equality def$(X,L)=n-2m$ implies that all inequalities $\dim(\Cal
D(X,V)) \leq \dim(p_1(I)) \leq \dim(I) = N-1-(2n-m)$ are
equalities so that $p_1$ is generically finite and $\Cal
D(X,V)=p_1(I)$. This completes the proof of (2.1.1).

If $n-2m=-1$ then $r=m$. Hence,  by (1.2.8), we get
$c_n(J_1(L))=2v(Y,E)-2+2c_m(E)$. Since $m \geq 3$, $v(Y,E) \geq
0$, see \cite{Fu, Thm. 4.1}, and $c_m(E)>0$ by \cite{BG, Cor.
1.2}. Hence $c_n(J_1(L))=0$ gives $v(Y,E)=0$ and $c_m(E)=1$. Then
(2.1.2) follows from \cite{Wi, Thm. 3.4}.\qed
\enddemo

\proclaim{Remark 2.2} {\rm Let us point out that the inclusion $\Cal
D(X,V) \subseteq p_1(I)$ provides a geometrical proof of the
inequality def$(X,L) \geq n-2m$. On the other hand, since $p_1(I)$
is irreducible, in the conditions of (2.1.1) the discriminant $\Cal
D(X,V)$ is irreducible. In particular, $\Cal J_n(X,V)=\emptyset$.}
\endproclaim

Let us observe that if $n-2m < -1$ then $\dim(I)=N-1-(n-2m)>N$.
Hence $p_1$ cannot be generically finite. In the following example
$n-2m=-1$ and $p_1(I) \subsetneq \Cal D(X,V)$, showing the
necessity of the hypothesis $n-2m \geq 0$ in (2.1.1).
\medskip

\proclaim{Example 2.3}{\rm Consider the Segre embedding of the
product $\Bbb P^1 \times \Bbb P^2$ in $\Bbb P^5$ as a scroll over
$\Bbb P^2$. Here $X=\Bbb P_{\Bbb P^2}(\Cal O_{\Bbb P^2}(1)^{\oplus
2})$ and the tautological line bundle $L$ is very ample. Then
$\dim(I)=5$ and $X^\vee =\Cal D(X,L) \subsetneq p_1(I)=(\Bbb
P^{5})^\vee$.}
\endproclaim

The hypothesis $r>1$ on the rank of $E$ in Proposition 2.1 is
crucial to get the inclusion $\Cal D(X,V) \subseteq p_1(I)$. The
following construction will allow us to produce examples (see
Example 2.4.3) in which $r=1$ and $\Cal D(X,V)$ is not contained
in $p_1(I)$, being in particular reducible. Let us also point out
that the following construction is of interest by its own. In fact
it relates the second order infinitesimal information of a variety
$Y \subset \Bbb P^N$ with the discriminant locus of a special
scroll over $Y$, namely its conormal variety, see \cite{R}.
\medskip

\noindent {\bf (2.4)} Consider an irreducible non-degenerate smooth
projective variety $Y \subset \Bbb P^{N}$ with $\dim(Y)=m$, such
that the twist of the normal bundle ${\Cal N}_{Y/\Bbb P^N}(-1)$ is
ample. For instance, this assumption is satisfied when $m=1$; or
when $m=2$ and Pic$(Y) \simeq \Bbb ZH$ ($H$ being the hyperplane
section); or when $Y \subset \Bbb P^N$ is a complete intersection.
As in \cite{LM2, Example 3.2} consider the conormal variety $X =\Bbb
P({\Cal N}_{Y/\Bbb P^N}(-1))\subset \Bbb P^N \times \Bbb P^{N \vee}$
(note that $n:=\dim(X)=N-1$) and the corresponding projections $\pi$
and $\pi_2$:
$$\matrix
\Bbb P^N \times \Bbb P^{N \vee} \supset &  \Bbb P({\Cal N}_{Y/\Bbb
P^N}(-1))= & X & \mapright{\pi} & Y \subset \Bbb P^N \cr & &
\mapdown{\pi_2} \cr & \phantom{................}\Bbb P^{N \vee}
\supset & \phantom{..} Y ^\vee
\endmatrix$$ where, by definition, $\pi_2(X)=Y^\vee$. The
following triplet is as in $(0.0)$:
$$(X,L= \pi_2^*\Cal O_{\Bbb P^{N\vee}}(1),V= \pi_2^*H^0(\Bbb P^{N
\vee},\Cal O_{\Bbb P^{N\vee}}(1))).$$

\proclaim{Remark 2.4.1} {\rm For a triplet as in (2.4) we have $\Cal
D_0(X,V)=(Y^\vee)^\vee=Y \subseteq \Cal D(X,V)$, so that $\text{\rm
def}_0(X,V)=N-1-\dim(\Cal D_0)=N-1-m=n-m.$ Moreover, by (2.1.1), if
$N-1\geq 2m$ then $\text{\rm def}(X,L) = n-2m$. Since $n-m>n-2m$, we
conclude that when $N \geq 2m+1$ there is a strict inclusion $\Cal
D_0(X,V) \subsetneq \Cal D(X,V)$.}
\endproclaim

The following lemma relates the jumping sets ${\Cal J}_i(X,V)$ of
a triplet as in $(2.4)$ with the hyperplane sections of $Y\subset
\Bbb P^N$ whose singularities are not general, i.e.,  not ordinary
quadratic of maximal rank. This shows the interaction between the
second order infinitesimal information of $Y \subset \Bbb P^N$ and
the discriminant locus of its conormal variety.

\proclaim{Lemma 2.4.2} For $(X,L,V)$ as in $(2.4)$ it follows that
$\Cal J_i=\{(y,H) \in X:$ {\rm the quadratic part of the defining
equation of $Y \cap H$ locally at $y$  has rank $\leq m-i\}$.
\endproclaim

\demo{Proof} Recall that, by definition, $\Cal J_i=\{(y,H) \in X:
\text{\rm rank}(d\pi_2(y,H))\leq \dim (X)-i\}$. Then the lemma is
a consequence of the local expression of $\pi_2$ at $(y,H)$, see
for example \cite{T, pp. 58--59} or \cite{R, Sect. 2}.
\qed\enddemo

Here is the promised example.

\proclaim{Example 2.4.3} {\rm Take $Y \subset \Bbb P^2$ a plane
curve and the triplet $(X,L,V)$ as in (2.4). Using the description
of $\Cal J_1$ of Lemma 2.4.2 we get that $\Cal J_1=\{(y,H) \in X:$
the quadratic part of the defining equation of $Y \cap H$ locally
at $y$ vanishes$\}$. Then, the description of the discriminant
locus of (0.2) is the following: $\Cal D(X,V)=\Cal D_0(X,V) \cup
\Cal D_1(X,V)$ where $\Cal D_0(X,V)=Y$ (see Remark 2.4.1) and
$\Cal D_1(X,V)$ is the union of all inflectional tangent lines to
$Y$. Any of these lines corresponds to the pencil of lines through
the cusp of $Y^{\vee}$ dualizing a flex of $Y$. In particular
$\Cal D(X,V)$ is reducible.}
\endproclaim

In connection with Example 2.4.3, let us include a remark about the
degree of the discriminant. Since in this example $X$ is isomorphic
to $Y$ then $L=(d-1)H$, where $d=$deg$(Y)$ and $H$ defines the
embedding $Y \subset \Bbb P^2$. Hence
$c_1(J_1(L))=c_1(J_1((d-1)H))=c_1(J_1(H))+2d(d-2)=(d-1)d+2d(d-2)$,
that is, $c_1(J_1(L))=d+3d(d-2).$ Since $\Cal D_0(X,V)=Y$ is an
irreducible component of $\Cal D(X,V)$ then, by \cite{LM3, Thm.
5.2}, the first summand $d$ corresponds to the degree of $\Cal
D_0(X,V)$. The second summand $3d(d-2)$ corresponds essentially to
the degree of $\Cal D_1(X,V)$. In fact, again by \cite{LM3, Thm.
5.2}, $3d(d-2)$ is the number of properly counted flexes of a degree
$d$ plane curve $Y \subset \Bbb P^2$. Concretely, each flex $f$ is
counted $m_f$ times, where $m_f$ is the Milnor number of the
singularity of the point $(f, T_{Y,f})$ in the element of $|V|$
through it. This is a natural generalization of the well known
Pl\"ucker formula, see \cite{K2, (I, 31)}.

As an application of Proposition 2.1 we can completely describe
the discriminant locus for a triplet as in (2.4) in the range
where 2.1 is meaningful.

\proclaim{Proposition 2.4.4} Let $(X,L,V)$ be a triplet as in
$(2.4)$. If $N \geq 2m+1$ then the discriminant locus $\Cal
D(X,V)$ can be identified with the tangent developable $TY \subset
\Bbb P^N$.
\endproclaim

\demo{Proof} We use the identification between $\Cal D(X,V)$ and
$p_1(I)$ of (2.1.1). For $y \in Y$ we have that
$\pi^{-1}(y)=\{{\Cal H} \in \Bbb P^{N\vee}: T_{Y,y} \subset {\Cal
H}\}$ by definition of conormal variety. Under the natural
identification between $\Bbb P^{N \vee \vee}$ and $\Bbb P^N=|V|$
the set $\{(H,y):\pi^{-1}(y) \subset H\} \subset |V| \times Y$ is
sent onto $T_{Y,y}$ by $p_1$. This shows that $TY=\Cal D(X,V)$.
\qed
\enddemo

\head 3 Characterization of def$(X,L)>n-4$
\endhead

In the classical context, varieties with big defect with respect
to their dimension tend to be scrolls. To be precise, if $L$ is
very ample and $\phi_V$ is an embedding, then def$(X,L) \geq
n-3>0$ if and only if $(X,L)$ is a scroll over a curve or
$(X,L)=(\Bbb P^n,{\Cal O}_{\Bbb P^n}(1))$, see \cite{LS, Cor.
3.4}, \cite{E, Thm. 3.2} and Landman's parity Theorem \cite{L},
\cite{K3, II (22)}. Moreover, when def$(X,L) \geq n-4>0$ we have
three possibilities, see \cite{M, Prop. 2.5}: either $(X,L)$ is a
scroll over a surface, or $X$ embedded by $|L|$ is $\Bbb G(1,4)
\subset \Bbb P^9$, i.e., the Pl\"ucker embedding of the
Grassmannian of lines in $\Bbb P^4$, or a smooth hyperplane
section of it. This shows that in this range examples different
from scrolls are very few.

We are going to use the computations of Section 1 to prove results
of this type in the ample and spanned case. In fact we also give new
proofs of the classical results when def$(X,L) \geq n-3$. To be
concrete: let $(X,L)$ be as in $(0.0)$, then def$(X,L)\leq n$ with
equality if and only if $(X,L)=(\Bbb P^n, {\Cal O}_{\Bbb P^n}(1))$,
\cite{LPS1, Thm. 2.8}. Furthermore, the next case def$(X,L) \ne n-1$
cannot occur, see \cite{LPS1, Thm. 2.8}. Here we deal with the
next-to-maximal cases def$(X,L)=n-i$, $i=2,3,4$. We characterize
scrolls over curves by the equality def$(X,L)=n-2$, we exclude the
possibility for the defect to be $n-3$ and we classify the pairs
$(X,L)$ in case def$(X,L)=n-4$ under the assumption that the Picard
number of $X$ is one.

The proofs in the classical setup use adjunction theory and either
topology or the linearity of the general contact locus. In our
proofs for the ample and spanned case we only use the adjunction
theoretic results on the ampleness of $K_X+(\dim(X)-i)L$ for
$i=0,1,2$ and $L$ ample and spanned together with the exact sequence
of \cite{LPS1, (2.8.2)}.

\proclaim{Lemma 3.0} Let $(X,L)$ be as in $(0.0)$ and suppose that
$(X,L)\ne (\Bbb P^n, \Cal O_{\Bbb P^n}(1))$. If {\rm def}$(X,L)>0$
then $K_X+(${\rm def}$(X,L)+2)L$ is not ample.
\endproclaim

\demo{Proof}  Suppose by contradiction that def$(X,L)=n-s>0$ with
$s$ an integer $n-1 \geq s \geq 0$ and $K_X+(n-s+2)L$ is ample. By
\cite{LM2, Lemma 1.13} and adjunction formula, just restricting
iteratively to general elements of $|L|$, we can suppose $n=s+1$,
def$(X,L)\geq 1$ and $K_X+3L$ ample. For any integer $r$, $1 \leq r
\leq n+1$ consider the filtration of $\wedge^r J_1(L)$ associated to
(0.3) (see \cite{H, pp. 127--128}): $\wedge^r J_1(L) = F_0 \supseteq
F_1 \supseteq \dots \supseteq F_r \supseteq F_{r+1}=0$, where $F_i /
F_{i+1} \cong \wedge^i (\Omega_X \otimes L) \otimes \wedge^{r-i}L$,
that is, $\Omega_X^{i}(rL)$ for $i = r-1,r$ and $0$ otherwise. This
provides the exact sequence:
$$ 0 \to \Omega_X^r(rL) \to \wedge^rJ_1(L) \to \Omega_X^{r-1}(rL)
\to 0.$$ Tensor this sequence by $M=-(K_X+ \alpha L)$ to get
$$0 \to \Omega_X^r(M+rL) \to \wedge^rJ_1(L)(M) \to \Omega_X^{r-1}(M+rL)
\to 0.$$ The line bundle $-(M+rL)=K_X+(\alpha-r)L$ is ample by
hypothesis if $\alpha-r \geq 3$. Then we can use the Nakano
vanishing theorem, \cite{La, Thm. 4.2.3, p. 250}, to get
$h^j(X,\Omega_X^r(M+rL))=0$ when $j+r<n$ and
$h^j(X,\Omega_X^{r-1}(M+rL))=0$ when $j+r-1<n$. This leads to this
vanishing: $$h^j(X, \wedge^rJ_1(L)(-(K_X+\alpha L)))=0,\; \alpha-r
\geq 3,\; j+r<n. \tag{3.0.1}$$ Since def$(X,L) \geq 1$ we can choose
a general line in $|L|$ not meeting $\Cal D(X,L)$. This allows us
construct the following exact sequence, exactly as in \cite{LPS1,
(2.8.2)}:
$$0 \to \Cal O_X^{\oplus 2} \to J_1(L) \to Q \to 0, \tag{3.0.2}$$
where $Q$ is a rank $n-1$ vector bundle on $X$. Now we claim that:
$$h^j(X, \wedge^{i}Q(-(K_X+(n+1)L)))=0,\; 0 \leq i \leq n-2, \; 0
\leq j \leq n-i-1. \tag{3.0.3}$$ Before proving the claim let us
show how it leads to a contradiction. Consider the exact sequence
dual of (0.3): $$0 \to \Cal O_X(-L) \to J_1(L)^\vee \to T_X(-L)
\to 0.$$ Since $(X,L)\ne (\Bbb P^n, \Cal O_{\Bbb P^n}(1))$ then
$h^0(X,T_X(-L))=0$, see \cite{W} or \cite{BS, Thm. 5.4.5}.
Moreover, $h^0(X,-L)=0$, whence $h^0(X,J_1(L)^\vee)=0$. Taking the
dual exact sequence of (3.0.2) we get $h^1(X, Q^\vee)>0$. Let us
recall that $Q^\vee =\wedge^{n-2}Q(-(K_X+(n+1)L)$ so that
$h^1(X,\wedge^{n-2}Q(-(K_X+(n+1)L))=h^1(X, Q^\vee)>0,$
contradicting the claim (3.0.3).

Now we can prove the claim by induction on $i$. For $i=0$, it
follows from the Kodaira vanishing theorem that $h^j(X,\Cal
O_{X}(-(K_X+(n+1)L))=0$ for $0 \leq j \leq n-1$, since $n+1 \geq 3$.

For $i=1$ the claim follows from (3.0.1), the Kodaira vanishing
theorem and the exact sequence obtained by twisting (3.0.2) by
$M=-(K_X+(n+1)L)$:
$$0 \to \Cal O_X^{\oplus 2}(M) \to J_1(L)(M) \to Q(M) \to 0.$$

For $i \geq 2$, using the filtration of $\wedge^i J_1(L)$
associated to (3.0.2), we get the following exact sequences:
$$0 \to F_1 \to \wedge^{i}J_1(L) \to \wedge^{i}Q \to 0,$$
$$0 \to \wedge^{i-2}Q \to F_1 \to (\wedge^{i-1}Q)^{\oplus 2} \to
0,$$ where $F_1$ is a vector bundle on $X$. Tensoring the second
exact sequence by $-(K_X+(n+1)L)$, by induction hypothesis, we get
$h^j(X,F_1(-(K_X+(n+1)L)))=0$ for $0 \leq j \leq n-i$. Tensoring the
first one by the same line bundle and using (3.0.1) we get the
claim.\qed \enddemo

\proclaim{Theorem 3.1} Let $(X,L)$ be as in $(0.0)$ with $n \geq
3$. The following are equivalent:

$(3.1.1)$ $X =\Bbb P_C(E)$ where $C$ is a smooth curve, $E$ is a
vector bundle on $C$ and $L=\Cal O_X(1)$ is the tautological
bundle;

$(3.1.2)$ {\rm def}$(X,L)=n-2$.
\endproclaim

\demo{Proof} If (3.1.1) holds then, since $n \geq 3$, def$(X,L)=n-2$
by (2.1.1).

Let us now prove that (3.1.2) implies (3.1.1) by induction on $n$.
Let us start with $n=3$, so that def$(X,L)=1$. Since by Lemma 3.0
$K_X+3L$ is not ample then we conclude by \cite{I, (1.3)}.

If $n>3$, consider a general element $X_1 \in |V|$, set
$L_1=L|_{X_1}$ and denote by $V_1$ the image of $V$ via the
restriction map $r:H^0(X,L) \to H^0(X_1,L_1)$. Then the triplet
$(X_1,L_1,V_1)$ is as in (0.0) and def$(X_1,L_1) \geq \dim(X_1)-2$
by [LM2, Lemma 1.13]. It cannot be $(X_1,L_1) =(\Bbb P^{n-1},\Cal
O_{\Bbb P^{n-1}}(1))$, otherwise it would be def$(X,L)=n$. Hence,
by induction, $X_1$ is a scroll over a smooth curve, $\dim(X_1)
\geq 3$ and thus we conclude by \cite{BS, Thm. 5.5.2}. \qed
\enddemo

\proclaim{Remarks 3.1.6} {\rm Let us point out that Theorem 3.1
improves \cite{LM2, Corollary 4.6}. Observe also that in (2.1.2) we
impose $m \geq 3$ to use the results of Fukuma. If $m=2$ then $n=3$
and $X$ is a scroll $\pi:X=\Bbb P_Y(E) \to Y$ over a smooth surface
$Y$, rank$(E)=2$. If def$(X,L)>0$ then def$(X,L)=\dim(X)-2$ so that,
by Theorem 3.1, $X$ must also be a scroll over a smooth curve $C$,
i.e., $\pi':X=\Bbb P_C(E') \to C$. These two structures of scroll on
$X$ lead easily to $Y=\Bbb P^2$ and $E={\Cal O}_{\Bbb
P^2}(1)^{\oplus 2}$. Then $(2.1.2)$ is also true when $m =2$.
Moreover, it is also true when $m=1$, even though the hypothesis
$r>1$ does not hold.}
\endproclaim

In the classical case Landman's parity theorem says that def$(X,L)$
and dim$(X)$ have the same parity, see reference above or \cite{T,
Thm. 7.4} and references inside. In particular def$(X,L)\ne n-3$. We
are going to prove this last result in the ample and spanned case.
Together with \cite{LPS1, Thm. 2.8} this seems to be a new evidence
towards Landman's parity theorem in the ample and spanned case.

Let us first recall the following definition, see \cite{I, (0.11)}
or \cite{BS}. Let $(X,L)$ be a polarized pair of dimension $n$,
i.e., $X$ is a smooth projective variety and $L$ is an ample line
bundle on $X$. An effective divisor $E \subset X$ is called a
$(-1)$-hyperplane if $E \simeq \Bbb P^{n-1}$, ${\Cal O}_X(L) \otimes
{\Cal O}_E={\Cal O}_E(1)$ and ${\Cal O}_X(E)\otimes {\Cal O}_E={\Cal
O}_E(-1)$. A smooth polarized variety $(X',L')$ is called a {\it
reduction of} $(X,L)$ if  is obtained contracting all the
$(-1)$-hyperplanes of $(X,L)$, that is, $\sigma: X \to X'$ is the
blowing up of a smooth variety $X'$ at $s$ distinct points $p_1,
\dots, p_s \in X'$, $E_i=\sigma^{-1}(p_i)$ and $L=\sigma^*L'-E$
where $E=E_1+\dots+E_s$ the total exceptional divisor. For $n \geq
3$ the contraction $\sigma$ is uniquely determined by $(X,L)$. Let
us point out that $|L'|$ can have base points. Anyway, the
discriminant locus $\cal D(X,L') \subseteq |L'|$ makes sense, though
in principle it could be that $\Cal D(X',L')=|L'|$. Let us write as
an abuse of notation def$(X',L')=\dim(|L'|)-1-\dim(\Cal D(X',L'))$.
In particular it could be that def$(X',L')=-1$. However this cannot
happen. In fact, we can prove the following lemma.

\proclaim{Lemma 3.2.0} Let $(X,L)$ be as in $(0.0)$, $n \geq 3$. We
get {\rm def}$(X,L)\leq$ {\rm def}$(X',L')$ for $(X',L')$ a
reduction of $(X,L)$.
\endproclaim

\demo{Proof} By definition $X$ is the blowing up $\sigma$ of $X'$ at
$s$ distinct points $p_1, \dots, p_s$ and $L=\sigma^*L'-E$ being
$E=E_1+\dots +E_s$ the total exceptional divisor. If $H \in |L|$ is
any smooth divisor then it does not contain any exceptional divisor
and its image $H'$ in $X'$ is also smooth. In fact
mult$_{p_i}(H')=$deg$(C_{p_i}H')=$deg$(H|_{E_i})=1$ for every $i$,
where $C_{p_i}H'$ denotes the tangent cone of $H'$ at $p_i$. Indeed
this implies that any linear subspace of $|L|$ containing only
smooth divisors maps to a linear subspace of $|L'|$ of the same
dimension containing only smooth divisors. Just by definition, see
(0.1), the defect cannot decrease and the lemma follows.\qed
\enddemo

\proclaim{Theorem 3.2} For $(X,L)$ as in $(0.0)$ with $n \geq 4$
we have {\rm def}$(X,L)\ne n-3$.
\endproclaim

\demo{Proof} Suppose by contradiction that def$(X,L)=n-3$. Let us
observe that for a general $X_1 \in |L|$, denoting $L_1 =L|_{X_1}$
we have def$(X_1,L_1) \geq \text{def}(X,L)-1$. By \cite{LPS1, Thm.
2.8} and Theorem 3.1 the strict inequality cannot hold. Hence, by
induction we can suppose $n=4$ and def$(X,L)=1$. Thus, by Theorem
3.1, $(X,L)$ is not a scroll over a smooth curve so that $K_X+4L$ is
ample. By Lemma 3.0 we know that $K_X+3L$ is not ample. Recall that
$L$ is ample and base point free. Hence one of the following holds,
see \cite{I, (1.6)}:

(3.2.2) $(X,L)$ is a Del Pezzo variety, i.e., $-K_X=3L$,

(3.2.3) $(X,L)$ is a hyperquadric fibration over a smooth curve,

(3.2.4) $(X,L)$ is a scroll over a smooth surface,

(3.2.5) there exists a reduction $(X',L')$ such that $K_{X'}+3L'$
is ample.

The case (3.2.2) can be easily excluded looking at the
classification of Del Pezzo manifolds, see for example \cite{F3,
Thm. 8.11, p. 72}. The case $L^4=1$ does not occur since $|L|$ is
base point free. If $L^4=2$ then $X$ is a double cover of $\Bbb P^4$
ramified along a smooth quartic so that ${\Cal D}_1$ is a
hypersurface. If $L^4 \geq 3$ then $L$ is very ample. Here, in order
to avoid the computations of $c_4(J_1(L))$, we can use the results
in the classical case to exclude this possibility.

The cases (3.2.3) and (3.2.4) come from the computation of the top
Chern class that we have done in Section 1. If $(X,L)$ is a
hyperquadric fibration then, by $(1.4)$, we get
$c_4(J_1(L))=(2e-b)+3(2e-5b)-4(1-g)>-4(1-g)$. In order to have
positive defect we need that $c_4(J_1(L))=0$. Hence $g=0$. Then
$c_4(J_1(L))=0$ implies $(2e-b)+3(2e-5b)=4$. Since $2e-b
>0$ and $2e-5 b \geq 0$, it is straightforward to check that
the vanishing of $c_4(J_1(L))$ cannot occur. If $(X,L)$ is a
scroll over a smooth surface then $c_4(J_1(L))=c_2(E)>0$ by
(1.2.2) and the ampleness of $E$, see \cite{K1, Thm. 3}.

Finally let us consider (3.2.5).  By Lemma 3.2.0 we can choose a
line $T' \subset |L'|$ such that any element in $T'$ is smooth.
Consider two distinct points in $T'$ and their corresponding
sections $s_0,s_1 \in |L'|$. They give rise to an exact sequence:
$$0 \to {\Cal O}_{X'}^{\oplus 2} \to J_1(L') \to Q' \to 0,$$ where
$Q'$ is a rank 3 vector bundle on $X'$. Let us show the
injectivity of the map $f: {\Cal O}_{X'}^{\oplus 2} \to J_1(L')$
in the vector space fibre. The arguments of \cite{LPS1, (2.6)}
work at points outside the base locus of $|L'|$. So we can confine
to the base points. Since Bs$|L'| \subseteq \{p_1, \dots, p_s\}$
we can suppose that Bs$|L'|=\{p_1, \dots, p_t \}$, where $t \leq
s$. Consider the point $p_1$ (the same arguments works for $p_i$,
$i \leq t$). Locally around $p_1$ the morphism $f: {\Cal
O}_{X'}^{\oplus 2} \to J_1(L')$ is given by
$$f_{p_1}(\lambda_0,\lambda_1)=\lambda_0(s_0(p_1),
ds_0(p_1))+\lambda_1(s_1(p_1),ds_1(p_1))=$$ $$=\lambda_0(0,
ds_0(p_1))+\lambda_1(0,ds_1(p_1))=(0,
\lambda_0ds_0(p_1)+\lambda_1ds_1(p_1)).$$ Since the line $T'$ does
not meet the vector subspace of sections of $L'$ singular at $p_1$
we get that $f_{p_1}$ is injective. Now we can reproduce the
arguments of Lemma 3.0 based on the non vanishing of
$h^1({Q'}^{\vee})$ to contradict the ampleness of $K_{X'}+3L'$.
\qed \enddemo

Now we study the next step def$(X,L)=n-4>0$. By Lemma 3.0 we get
$K_X+(n-2)L$ is not ample. In the case def$(X,L)=n- 3$, when
$K_X+(n-1)L$ is not ample we know that $K_{X'}+(n-1)L'$ is ample
on the reduction $(X',L')$ and we can do a complete classification
using Lemma 3.2.0. For the case def$(X,L)=n- 4$ we cannot
reproduce the same argument because if $K_X+(n-2)L$ is not ample
then we can only get the semi-ampleness (that is, global
generation of a power) of the adjoint bundle $K_{X'}+(n-2)L'$ on
the reduction. However, we can prove the following under the
assumption that the Picard number is one.

\proclaim{Proposition 3.3} Let $(X,L)$ be as in $(0.0)$ with {\rm
def}$(X,L)=n-4>0$. If the Picard number of $X$ is equal to one
then either $X=\Bbb G(1,4)$ is the grassmannian of lines in $\Bbb
P^4$ and $L$ defines the Pl\"ucker embedding in $\Bbb P^9$ or $X$
is a smooth hyperplane section $Y$ of $\Bbb G(1,4) \subset \Bbb
P^9$ and $L=\Cal O_{\Bbb P^9}(1) \otimes \Cal O_Y$.
\endproclaim

\demo{Proof} Recall that $K_X+(n-2)L$ is not ample by Lemma 3.0.
Since Pic$(X)=\Bbb Z H$ with $H$ an ample line bundle, we easily see
that $H=L$. Moreover, either $(X,L)$ is a Del Pezzo manifold or
$K_{X}+(n-2)L$ is semi-ample, $(X,L)$ coinciding with its reduction,
by \cite{I, (1.7)}. In the former case we can use the classification
of Del Pezzo manifolds, \cite{F3, Thm. 8.11, p. 72}, to conclude
that $(X,L)$ is as in the statement. In the latter case
$-K_X=(n-2)L$, that is, $(X,L)$ is a {\it Mukai variety}. Recall
that $|L|$ is base point free by hypothesis. As a consequence of the
study of the anticanonical linear system for Fano threefolds of
\cite{Is} we get the following, see \cite{Mu, Prop. 1}: either $L$
is very ample, or $|L|$ defines a degree two morphism onto $\Bbb
P^n$, or onto a quadric $\Bbb Q^n \subset \Bbb P^{n+1}$. By \cite{M,
Prop. 2.5} $L$ cannot be very ample. If we have a degree two map
onto $\Bbb P^n$, $\phi_L:X \to \Bbb P^n$, then the branch locus of
$\phi_L$ is a smooth hypersurface (of degree $>1$) and its dual
variety is contained in ${\Cal D}(X,L)$. This contradicts
def$(X,L)>0$. If we have a degree two map $\phi_L:X \to \Bbb Q^n
\subset \Bbb P^{n+1}$ onto a quadric then $(X,L)$ has $\Delta$-genus
2 and  fits into \cite{F3, (10.8.2), p. 89}, so that $\Bbb Q^n
\subset \Bbb P^{n+1}$ is smooth, being $n \geq 3$. But then $\Bbb
Q^{n \vee} \simeq \Bbb Q^n \subseteq \Cal D(X,L)$ contradicting
def$(X,L)>0$. \qed
\enddemo

The following table summarizes the results on the classification
of varieties with high defect with respect to their dimension in
comparison with the classical case.
$$\matrix \noalign{\hrule \smallskip}  {\text {\rm def}}(X,L)>0 & n & n-1 & n-2 & n-3 & n-4
\cr \noalign{\smallskip\hrule \smallskip} L \text{ ample and} &
(\Bbb P^n, {\Cal O}_{\Bbb P^n}(1)) & {\text{Not}} & \text{Scroll
over a} &\text{Not } & \Bbb G(1,4) \subset \Bbb P^9\cr \text{spanned
by } V & &\text{possible} & \text{ curve} & \text{possible} & \Bbb
G(1,4) \cap H \cr  & & & &  & \text{if Pic}(X)=\Bbb Z  \cr
\noalign{\smallskip\hrule
\smallskip}\phi_V \text{ an} & (\Bbb P^n, {\Cal O}_{\Bbb P^n}(1)) &
{\text{Not}} & \text{Scroll over a} &\text{Not} & \text{Scroll
over a} \cr \text{embedding} & &\text{possible} & \text{curve} &
\text{possible} & \text{surface}\cr & &
 & & & \Bbb G(1,4) \subset \Bbb P^9\cr &
& & & & \Bbb G(1,4) \cap H \cr \noalign{\smallskip\hrule
\smallskip}
\endmatrix $$

\enddocument